\documentclass[reqno,12pt]{amsart}
\usepackage{amsthm}
\usepackage{amsmath,mathrsfs}
\usepackage{color}
\usepackage{amssymb}
\usepackage{hyperref}
\usepackage{enumerate}
\usepackage{latexsym,array}
\usepackage{amsfonts}
\usepackage{shadow}

\newtheorem{Pa}{Paper}[section]
\newtheorem{Tm}[Pa]{{\bf Theorem}}
\newtheorem{La}[Pa]{{\bf Lemma}}
\newtheorem{Dn}[Pa]{{\bf Definition}}
\newtheorem{Cy}[Pa]{{\bf Corollary}}

\newtheorem{Rk}[Pa]{{\bf Remark}}
\newtheorem{Pn}[Pa]{{\bf Proposition}}
\newtheorem{Pb}[Pa]{{\bf Problem}}

\def\hh{\mathbb{H}}

\author[K. Abu-Ghanem]{Khaled Abu-Ghanem}
\address{(KA) Department of Mathematics\\
Ben-Gurion University of the Negev\\
Beer-Sheva 84105 Israel} \email{khaledab@post.bgu.ac.il}
\author[D. Alpay]{Daniel Alpay}
\address{(DA) Department of Mathematics\\
Ben-Gurion University of the Negev\\
Beer-Sheva 84105 Israel} \email{dany@math.bgu.ac.il}

\author[F. Colombo]{Fabrizio Colombo}
\address{(FC) Politecnico di
Milano\\Dipartimento di Matematica\\Via E. Bonardi, 9\\20133
Milano, Italy}
\email{fabrizio.colombo@polimi.it}

\author[D. P. Kimsey]{David P. Kimsey}
\address{(DPK) Department of Mathematics\\
Ben-Gurion University of the Negev\\
Beer-Sheva 84105 Israel}
\email{dpkimsey@gmail.com}

\author[I. Sabadini]{Irene Sabadini}
\address{(IS) Politecnico di
Milano\\Dipartimento di Matematica\\Via E. Bonardi, 9\\20133
Milano, Italy}
\email{irene.sabadini@polimi.it}
\title[Boundary interpolation in the quaternionic setting]
{Boundary interpolation for slice hyperholomorphic Schur functions} \oddsidemargin
0.2in \evensidemargin 0.2in \topmargin -0.5in \textwidth
15.5truecm \textheight 23truecm

\def\R{\mathbb R}

\def\(s){\mathscr S(\R\times\R)}



\thanks{D. Alpay thanks the Earl Katz family for endowing the chair
which supported his research, and the Binational Science
Foundation Grant number 2010117. F. Colombo and I. Sabadini
acknowledge the Center for Advanced Studies of the Mathematical
Department of the Ben-Gurion University of the Negev for the
support and the kind hospitality during the period in which part
of this paper has been written. D. P. Kimsey acknowledges funding from a
Kreitman postdoctoral fellowship.}
\begin{document}
\parindent 0cm
\begin{abstract}
A boundary Nevanlinna-Pick interpolation problem is posed and solved in the
quaternionic setting. Given nonnegative real numbers $\kappa_1, \ldots, \kappa_N$,
quaternions \\ $p_1, \ldots, p_N$ all of modulus $1$, so that the $2$-spheres
determined by each point do not intersect and $p_u \neq 1$ for $u = 1,\ldots, N$,
and quaternions $s_1, \ldots, s_N$, we wish to find a slice hyperholomorphic Schur
function $s$ so that
$$\lim_{\substack{r\rightarrow 1\\ r\in(0,1)}} s(r p_u) = s_u\quad {\rm for} \quad u=1,\ldots, N,$$
and
$$\lim_{\substack{r\rightarrow 1\\ r\in(0,1)}}\frac{1-s(rp_u)\overline{s_u}}{1-r}\le\kappa_u,\quad
{\rm for} \quad u=1,\ldots, N.$$
Our arguments relies on the
theory of slice hyperholomorphic functions and reproducing kernel
Hilbert spaces.
\end{abstract}
\maketitle \tableofcontents

\parindent 0cm

\section{Introduction}
\setcounter{equation}{0}

In the paper \cite{2013arXiv1308.2658A} the Nevanlinna-Pick interpolation problem for
slice hyperholomorphic Schur functions has been solved using the
FMI (fundamental matrix inequality) method (see \cite{kky} for details). By a Schur function we mean a function $f$ which is slice hyperholomorphic on the open unit ball $\mathbb B_1$ of the quaternions and is bounded in modulus by $1$, i.e. $\sup_{p\in\mathbb B_1}|f(p)|\leq 1$. In the present paper
we solve a boundary interpolation problem for slice hyperholomorphic functions using the reproducing kernel Hilbert space method based on de Branges-Rovnyak spaces.
We refer the reader to \cite{abds2,abds3,Dym_CBMS} for more information on the reproducing kernel
Hilbert space approach to interpolation problems.\\

We state the problem we will solve in this paper and introduce some notation and definitions. Let us denote by  $\mathbb B_1$ and $\mathbb H_1$, the open unit ball and the unit sphere of $\mathbb H$, respectively.
For a given element $p\in\mathbb H$ we denote by $[p]$ the associated 2-sphere:
\[
[p]=\left\{qpq^{-1}: q\in\mathbb H\setminus\left\{0\right\}\right\}.
\]
Recall that two quaternions belong to the same sphere if and only if they have the same modulus and the same real part.

\begin{Pb}
\label{pb1}
Given $p_1,\ldots,p_N\in\mathbb H_1\setminus \left\{1\right\}$ such that $[p_u]\cap[p_v]=\emptyset$ for $u\not = v$ (the interpolation nodes),  $s_1,\ldots,
s_N\in\mathbb H_1$, and
$\kappa_1,\ldots, \kappa_N\in [0,\infty)$, find a necessary and
sufficient condition for a slice hyperholomorphic Schur function
$s$ to exist such that the conditions
\begin{eqnarray}
\label{inter1}
\lim_{\substack{r\rightarrow 1\\ r\in(0,1)}} s(rp_u)&=&s_u,\\
\lim_{\substack{r\rightarrow 1\\
r\in(0,1)}}\frac{1-s(rp_u)\overline{s_u}}{1-r}&\le& \kappa_u
\label{inter2}
\end{eqnarray}
hold for $u=1, \ldots N$, and describe the set of all Schur functions satisfying \eqref{inter1}-\eqref{inter2} when this condition is in force.
\end{Pb}

We note that \eqref{inter1}-\eqref{inter2} imply that
\begin{equation}
\lim_{\substack{r\rightarrow 1\\
r\in(0,1)}}\frac{1-|s(rp_u)|^2}{1-r^2}\le \kappa_u,\quad
u=1,\ldots, N, \label{wertyu} \end{equation} since
\begin{equation}
\label{richelieu-drouot1}
\frac{1-|s(rp_u)|^2}{1-r^2}=\frac{1-s(rp_u)\overline{s_u}}{(1-r)(1+r)}+
(s(rp_u)\overline{s_u})\frac{1-s_u\overline{s(rp_u)}}{(1-r)(1+r)}.
\end{equation}

We also note that the fact that the limits \eqref{wertyu} is part
of the requirement in the interpolation problem (in the complex case, the
corresponding limit is well-known to be non-negative).\\

As it appears from the statement of Problem \ref{pb1}, there is a
major difference with the complex case. Here we have to require
that not only the interpolation points are distinct, but also the
spheres they determine. The fact that this hypothesis is
necessary, and cannot be avoided, can be intuitively justified by
the fact that the $S$-spectrum of a matrix, or in general of an
operator (see Definition \ref{defspscandres}), consists of
spheres (which may reduce to real points). It is important to
note that the notion of $S$-spectrum of a matrix $T$ coincides
with the set of right eigenvalues of $T$, i.e. the set of
$\lambda \in \mathbb{H}$ so that $Tx = x \lambda$ for a nonzero
vector $x$. \\

Another major difference is the lack of a Carath\'eodory theorem
(see e.g. \cite[p. 48]{sarason94}) in the quaternionic setting.\\

Part of the arguments follow the classical case, taking into
account the noncommutativity of the quaternions. As we shall see,
even though the structure of the proof follows the the arguments
from \cite{adubi1}, it is necessary to suitably adapt the
arugment to the quaternionic setting and often the needed
modifications are not immediate.\\

The paper consists of five sections, besides the introduction. In
Section 2, we recall some basic material on slice
hyperholomorphic functions which will be needed in the sequel.
Section 3 illustrates the strategy and the various steps we will
follow to solve Problem \ref{pb1}. Section 4 contains detailed
proofs of these steps and Section 5 deals with the degenerate
case. Section 6 deals with an analogue of Carath\'eodory's
theorem in the quaternionic setting.

\section{Some preliminaries}
\setcounter{equation}{0}
In this section we collect some basic results, which will be used in the sequel.
Let $\hh$ be the real associative algebra of quaternions
with respect to the basis $\{1, i,j,k \}$
satisfying the relations
$
i^2=j^2=k^2=-1,\
 ij =-ji =k,\
jk =-kj =i ,
 \  ki =-ik =j .
$
A quaternion $p$ is denoted by $p=x_0+ix_1+jx_2+kx_3$,
$x_\ell\in \mathbb{R}$, $\ell=0,\ldots, 3$, its conjugate is
$\bar p=x_0-ix_1-jx_2-kx_3$, and the norm of a quaternion is such that $|p|^2=p\overline{p}$.
A quaternion $p$ can be written as $p={\rm Re}(p)+\underline{p}$ where the real part ${\rm Re}(p)$ is $x_0$  and $\underline{p} = i x_1 + j x_2 + k x_3$.
The symbol $\mathbb{S}$ denotes the 2-sphere of purely imaginary unit quaternions, i.e.
$$
\mathbb{S}=\{ \underline{p}=ix_1+jx_2+kx_3\ |\  x_1^2+x_2^2+x_3^2=1\}.
$$
Note that if $I\in\mathbb S$ then $I^2=-1$. Any nonreal quaternion $p=x_0+ix_1+jx_2+kx_3$ uniquely determines an element $I_p=(ix_1+jx_2+kx_3)/|ix_1+jx_2+kx_3|\in\mathbb S$. If $p=x_0\in\mathbb R$ then $p=x_0+I0$ for all $I\in\mathbb S$. Given $p\in\mathbb H$ we can write $p=p_0+I_pp_1$ and the 2-sphere $[p]$ coincides
 with the set of all elements of the form $p_0+Jp_1$ when $J$ varies in $\mathbb{S}$.
The set $[p]$ is reduces to the point $p$ if and only if $p\in\mathbb{R}$.\\
We now recall the definition of a slice hyperholomorphic function, for more details see \cite{MR2752913}.
\begin{Dn}
{\rm
Let $\Omega\subseteq\hh$ be an open set and let
$f:\ \Omega\to\hh$ be a real differentiable function. Let
$I\in\mathbb{S}$ and let $f_I$ be the restriction of $f$ to the
complex plane $\mathbb{C}_I := \mathbb{R}+I\mathbb{R}$ passing through $1$
and $I$ and denote by $x+Iy$ an element on $\mathbb{C}_I$.
 We say that $f$ is a left slice hyperholomorphic (or slice hyperholomorphic, for short) function
  in $\Omega$ if, for every
$I\in\mathbb{S}$, we have
$$
\frac{1}{2}\left(\frac{\partial }{\partial x}+I\frac{\partial
}{\partial y}\right)f_I(x+Iy)=0.
$$
 We say that $f$ is a right slice hyperholomorphic function
 in $\Omega$ if, for every
$I\in\mathbb{S}$, we have
$$
\frac{1}{2}\left(\frac{\partial }{\partial x}f_I(x+Iy)+\frac{\partial
}{\partial y}f_I(x+Iy) I\right)=0.
$$
}
\end{Dn}
Slice hyperholomorphic functions have nice properties on some particular open sets which are defined below.
\begin{Dn}
{\rm
Let $\Omega$ be a domain in $\mathbb{H}$.
We say that $\Omega$ is a
slice domain (s-domain for short) if $\Omega \cap \mathbb{R}$ is non empty and if
$\Omega\cap \mathbb{C}_I$ is a domain in $\mathbb{C}_I$ for all $I \in \mathbb{S}$.
We say that $\Omega$ is
axially symmetric if, for all $p \in \Omega$, the
sphere $[p]$ is contained in $\Omega$.
}
\end{Dn}
On an axially symmetric s-domain $\Omega$, a slice hyperholomorphic function satisfies the following
 formula, which is called the Structure formula or the Representation formula (see \cite[Theorem 4.3.2]{MR2752913}):
\begin{equation}\label{repr}
f(x+J y)=\frac 12 \left[ f(x+Iy) +f(x-Iy) + J I
(f(x-Iy)- f(x+Iy))\right].
\end{equation}
Formula \eqref{repr} is useful as it allows one to extend a holomorphic map
  $h:\ \Omega\subseteq \mathbb{C}\cong \mathbb{C}_I\to \mathbb H$ to a slice hyperholomorphic function. Let $U_\Omega$ be the axially symmetric completion of $\Omega$, i.e.
  $$U_\Omega=\bigcup_{J\in\mathbb S, \, x+Iy\in\Omega} \{x+Jy\}.$$
The left slice hyperholomorphic extension ${\rm ext}(h): \ U_\Omega\subseteq \mathbb H \to\mathbb H$ of $h$ is the function defined as (see \cite{MR2752913}):
\begin{equation}\label{ext}
{\rm ext}(h)(x+J y)=\frac 12 \left[ h(x+Iy) +h(x-Iy) + J I
(h(x-Iy)-h(x+Iy))\right].
\end{equation}
It is immediate that  ${\rm ext}(h+g)={\rm ext}(h)+{\rm ext}(g)$
and that if $h(z)=\sum_{n=0}^\infty h_n(z)$ then ${\rm
ext}(h)(z)=\sum_{n=0}^\infty {\rm ext}(h_n)(z)$.

Two left (resp. right) slice hyperholomorphic functions can be multiplied, on an axially symmetric s-domain, using the so called $\star$-product (resp. $\star_r$-product) in order to obtain another left (resp. right) slice hyperholomorphic function.
\\
Let $f,g:\ \Omega \subseteq\mathbb{H}$ be slice hyperholomorphic functions. Their restrictions to the complex plane $\mathbb{C}_I$ can be written as
$f_I(z)=F(z)+G(z)J$,
$g_I(z)=H(z)+L(z)J$ where $J\in\mathbb{S}$, $J\perp I$, i.e. $IJ = -JI$. The functions $F$,
$G$, $H$, $L$ are holomorphic functions of the variable $z\in
\Omega \cap \mathbb{C}_I$, see \cite[p. 117]{MR2752913}.
We have the following:
\begin{Dn}
{\rm
Let $f$ and $g$ be slice hyperholomorphic functions defined on an axially symmetric s-domain $\Omega\subseteq\mathbb{H}$.
The $\star$-product of  $f$ and $g$ is defined as the unique
left slice hyperholomorphic function on $\Omega$ whose restriction to the
complex plane $\mathbb{C}_I$ is given by
\begin{equation}\label{starproduct}
\begin{split}
(f\star g)_I(z)&=(F(z)+G(z)J)\star(H(z)+L(z)J)\\
&=
(F(z)H(z)-G(z)\overline{L(\bar z)})+(G(z)\overline{H(\bar z)}+F(z)L(z))J.
\end{split}
\end{equation}
}
\end{Dn}
 If $f$ and $g$ are slice hyperholomorphic on a ball with center at the origin, they can be expressed in a power series, i.e. $f(p)=\sum_{n=0}^\infty p^n a_n$ and $g(p)=\sum_{n=0}^\infty p^n b_n$. Thus $(f\star g)(p)=\sum_{n=0}^\infty p^n c_n$, where
$c_n=\sum_{r=0}^na_rb_{n-r}$ is obtained by convolution on the coefficients.
For the construction of the  $\star$-product of right slice hyperholomorphic functions and for more information on the $\star$-product, we refer the reader to \cite{MR3127378, MR2752913}.
\\ Given a slice hyperholomorphic function, it is possible to define its slice hyperholomorphic reciprocal, see  \cite{MR2752913}. Here we limit ourselves to the case in which $f$ admits the power series expansion $f(p)=\sum_{n=0}^\infty p^n a_n$. In this case we set
$$
f^c(p)=\sum_{n=0}^\infty p^n \bar a_n,\qquad  f^s(p)=(f^c\star f)(p
)=\sum_{n=0}^\infty p^nc_n,\quad
c_n=\sum_{r=0}^n a_r\bar a_{n-r},
$$
so that
the left slice hyperholomorphic reciprocal of $f$
is defined as
$$
f^{-\star}:=(f^s)^{-1}f^c.
$$
In the general case, this formula is still valid with $f^s$, $f^c$ suitably defined.
\begin{Rk}{\rm  Let $k(p,q)$ be a function left slice hyperholomorphic in $p$ and right slice hyperholomorphic in $\bar q$. When taking the $\star$-product of a  function $f(p)$ slice hyperholomorphic in the variable $p$ with a function $k(p,q)$, we will write $f(p)\star k(p,q)$ meaning that the $\star$-product is taken with respect to the variable $p$;  similarly, the $\star_r$-product of $k(p,q)$ with  functions right slice hyperholomorphic in the variable $\bar q$ is always taken with respect to $\bar q$.
}
\end{Rk}

The following proposition is taken from \cite[Proposition 4.3]{MR3127378}, where a proof can be found.

\begin{Pn}
Let $\mathcal H(K_1)$ and $\mathcal H(K_2)$ be two reproducing kernel Hilbert spaces of
$\mathbb H^m$ and $\mathbb H^n$-valued slice hyperholomorphic functions in $\Omega$, with reproducing kernels
$K_1$ and $K_2$, respectively. Let $R$ be a $\mathbb H^{n\times m}$-valued function slice-hyperholomorphic in $\Omega$. Then the operator
of left $\star$-multiplication
\[
M_R\,:\, \,\, f\,\,\mapsto\,\, R\star f
\]
is continuous from  $\mathcal H(K_1)$ into $\mathcal H(K_2)$ if and only if the kernel
\[
K_2(p,q)-R(p)\star K_1(q, p)\star_rR(q)^*
\]
is positive definite in $\Omega$. Furthermore
\begin{equation}
\label{grenelle}
M_R^*(K_2(\cdot, q)d)=K_1(\cdot, q)\star_r R(q)^*d,\quad d\in\mathbb H^n.
\end{equation}
\label{la-seine}
\end{Pn}
Let us recall a few facts on the $S$-spectrum and on the $S$-resolvent operator.
\begin{Dn}\label{defspscandres}
{\rm
Let $A$ be a bounded quaternionic linear operator acting on a quaternionic, two sided, Banach space $V$.
We define the $S$-spectrum $\sigma_S(A)$ of $A$  as:
$$
\sigma_S(A)=\{ \text{$s\in \mathbb{H}$ : $A^2-2 {\rm Re}\,(s) A+|s|^2\mathcal{I}$ is not invertible}  \},
$$
where $\mathcal I$ denotes the identity operator on $V$.
The $S$-resolvent set $\rho_S(A)$ is defined as
$\rho_S(A)=\mathbb{H}\setminus\sigma_S(A)$.
}
\end{Dn}
From Definition \ref{defspscandres} it follows that the $S$-spectrum consists of spheres (which may reduce to real points).\\
The definition of $S$-spectrum arises from the following:
\begin{Pn}\label{Ssinistro}
Let $A$ be a bounded quaternionic linear operator acting on a quaternionic, two sided, Banach space $V$.
Then,  for $\|A\|< |p|$,  we have
\begin{equation}\label{SresolvR}
\sum_{n= 0}^\infty  s^{-1-n}A^n=-
(A-   \overline{s}\mathcal{I})(A^2-2{\rm Re}(s) A+|s|^2 \mathcal{I}  )^{-1}.
\end{equation}
\end{Pn}
\begin{Dn}
{\rm
The operator
\begin{equation}\label{SresolvoperatorRdd}
S_R^{-1}(s,A):=-(A-   \overline{s}\mathcal{I})(A^2-2{\rm Re}(s) A+|s|^2 \mathcal{I}  )^{-1},
\end{equation}
is called the right $S$-resolvent operator.
}
\end{Dn}
The right $S$-resolvent operator is obviously defined for $s\in\rho_S(A)$.\\
In the sequel we will be in need of the result below:
\begin{Pn}
Let $V$ be a two sided quaternionic Banach space and let
$A$ be a bounded right linear operator from $V$ into
itself. Then, for $|p| \,\|A\|< 1$ we have
\begin{equation}
\label{eq:oberkampf_ligne_5}
\sum_{n=0}^\infty p^n A^n =( \mathcal{I}  -\bar p A)(|p|^2A^2-2 {\rm Re}(p) A+ \mathcal{I}  )^{-1}.
\end{equation}
\end{Pn}
Another way to write the operator on the right hand side of (\ref{eq:oberkampf_ligne_5}) is to observe that it corresponds to the function one obtains by constructing the right $\star$-reciprocal of the function $f(q)=(1-pq)$. Upon computing $f^{-\star}(A)$ using the quaternionic functional calculus, see \cite{MR2752913}, one can write:
\begin{equation}\label{eq:oberkampf_ligne_55}
(\mathcal{I}-pA)^{-\star}=\sum_{n=0}^\infty p^n A^n.
\end{equation}

Finally, we mention a result which is a restatement of  \cite[Proposition 2.22]{2013arXiv1310.1035A} and which contains an identity that will be crucial in the sequel.
\begin{Pn}
Let $p\in\mathbb H$, $1/p\in\rho_S(A)$ and $(G,A)\in\mathbb H^{n\times m}\times{\mathbb H}^{m\times m}$. Then
\begin{equation}
\label{magic}
\sum_{t=0}^\infty p^tGA^t=(G-\overline{p}GA)(\mathcal{I}_m-2{\rm Re}(p)A+|p|^2A^2)^{-1},
\end{equation}
where $\mathcal{I}_m$ denotes the $m \times m$ identity matrix.
\end{Pn}
\begin{Rk}{\rm
We note that if $m=1$ then $A$ is a quaternion $a$ and the condition $1/p\in\rho_S(A)$ translates to the condition $1/p\not\in[a]$.
}
\end{Rk}

\section{The main result and the strategy}
\setcounter{equation}{0}
For the convenience of the reader we recall the main steps of the reproducing kernel method. We first introduce some notation. We set
\begin{equation}\label{Dn:ACJ}
A={\rm diag}\,(\overline{p_1},\ldots,\overline{p_N})\in\mathbb H^{N\times N},\quad C=\begin{pmatrix}1&\cdots &1\\
\overline{s_1}&\cdots &\overline{s_N}\end{pmatrix}\in\mathbb H^{2\times N},
\end{equation}
and
\[
\mathcal{J}=\begin{pmatrix}1&0\\0&-1\end{pmatrix}\in \mathbb R^{2\times 2}.
\]
Consider the matrix equation
\begin{equation}
\label{steineq}
P-A^*PA=C^*\mathcal{J}C
\end{equation}
where the unknown is $P\in\mathbb H^{N\times N}$.
The off diagonal entries of the matrix equation
are uniquely determined by the
equation
\begin{equation}
\label{puv}
P_{uv}-p_uP_{uv}\overline{p_v}=1-s_u\overline{s_v}
\end{equation}
but, in view of the following lemma the diagonal entries can be arbitrary:
\begin{La}\label{La:ph-hq}
Let $p$ and $q$ be quaternions of modulus $1$. Then, the equation
\begin{equation}
ph-hq=0,
\label{trivial}
\end{equation}
where $h\in\mathbb H$,
has the only solution $h=0$ if and only if ${\rm Re}(p)\not={\rm Re}(q)$, that is, if and only if $[p]\cap[q]=\emptyset$.
\end{La}

\begin{proof} If \eqref{trivial} has a solution $h\not=0$, then $p=hqh^{-1}$ and so $p$ and $q$ are in the same sphere. So a necessary
condition for \eqref{trivial} to have only $h=0$ as solution is that  $[p]\cap[q]=\emptyset$. We now show that this condition is also
sufficient.  Let $p=z_1+z_2j$ and $q=w_1+w_2j$, where $z_1,z_2,w_1,w_2\in\mathbb C$.
Since ${\rm Re}(p)\not={\rm Re}(q)$ we have
\begin{equation}
\label{complicated}
{\rm Re}( z_1)\pm i\sqrt{1-({\rm Re}( z_1))^2}\not={\rm Re}( w_1)\pm i\sqrt{1-({\rm Re}( w_1))^2}.
\end{equation}
We now introduce the
injective ring homomorphism $\chi : \mathbb H \to \mathbb C^{2 \times 2}$ given by
\begin{equation}
\label{eq:Oct27jkl1}
\chi(p) =  \begin{pmatrix} z_1 & z_2 \\ - \overline{z}_2 & \overline{z}_1 \end{pmatrix}.
\end{equation}
Using the map $\chi$, equation \eqref{trivial} becomes
\begin{equation}
\label{detour}
\chi(p)\chi(h)-\chi(h)\chi(q)=0.
\end{equation}
The eigenvalues of $\chi(p)$ are the solutions of
\[
\lambda^2-2({\rm Re}( z_1))\lambda+1=0,
\]
that is $\lambda={\rm Re}( z_1)\pm i\sqrt{1-({\rm Re}( z_1))^2}$, and similarly for $\chi(q)$.
By a well known result on matrix equations (see e.g., Corollary 4.4.7 in \cite{HornJohnson}), equation (\ref{detour}) has only the solution $\chi(h)=0$ if and only if $\lambda-\mu\not=0$ for all possible choices of eigenvalues of $\chi(p)$ and $\chi(q)$, and this condition holds in view of \eqref{complicated}. So the only solution
of \eqref{detour} is $h=0$.
\end{proof}
We denote by $P$ the $N\times N$ Hermitian matrix with entries $P_{uv}$ given by
\eqref{puv} for $u\not= v$ and with diagonal entries equal to $P_{uu}=\kappa_u$, $u,v=1,\ldots, N$. When $P$ is invertible we define
\begin{equation}
\label{parislelouvre}
\Theta(p)=\mathcal{I}_2-(1-p)\star C\star(\mathcal{I}_N-pA)^{-\star}P^{-1}(\mathcal{I}_N-A)^{-*}C^*\mathcal{J}=\begin{pmatrix}a(p)&b(p)\\ c(p)&d(p)\end{pmatrix}.
\end{equation}
Note that $\Theta$ is well defined in $\mathbb B_1$ since we assumed that the interpolation nodes $p_u$ are all different from $1$.
Finally we denote by $\mathcal M$ the span of the columns of the function
\begin{equation}
\label{Fp}
F(p)=C\star (\mathcal{I}_N-pA)^{-\star}
=\sum_{t=0}^\infty p^tCA^t,
\end{equation}
and endow $\mathcal M$ with the Hermitian form
\[
[F(p)c,F(p)d]_{\mathcal M}=d^*Pc,\quad c,d\in\mathbb H^N.
\]
We prove the following theorem.

\begin{Tm}\mbox{}\\
$(1)$ There always exists  a Schur function so that \eqref{inter1}
holds.\\
$(2)$ Fix $\kappa_1,\ldots, \kappa_N \geq 0$ and assume $P> 0$. Any
solution of Problem \ref{pb1} is of the form
\begin{equation}
\label{lft1}
s(p)=(a(p)\star e(p)+b(p))\star(c(p)\star
e(p)+d(p))^{-\star},
\end{equation}
where $a,b,c,d$ are as in (\ref{parislelouvre}) and $e$ is a
slice hyperholomorphic Schur function.\\
$(3)$ Conversely, any function of the form \eqref{lft1} satisfies
\eqref{inter1}. If
\begin{equation}
\label{limits}
\lim_{\substack{r\rightarrow 1\\
r\in(0,1)}}\frac{1-s(rp_u)\overline{s_u}}{1-r}
\end{equation}
exists and is real, then $s$ satisfies \eqref{inter2}.\\
$(4)$ If $e$ is a unitary constant, then the limit
\eqref{limits} exists (but are not necessarily real) and satisfies
\begin{equation}
\label{bastille}
\frac{|\beta_u-\overline{p_u}\beta_u\overline{p_u}|^2}{|1-\overline{p_u}^2|}\le
({\rm Re}\,\beta_u)\kappa_u.
\end{equation}
\label{tm1}
\end{Tm}

The strategy of the proof is as follows:\\

STEP 1: {\sl The condition $P\ge 0$ is necessary for Problem \ref{pb1} to have a solution.}\\

STEP 2: {\sl Assume that $s$ is a solution of Problem \ref{pb1}. Then  the map $M_{\begin{pmatrix}1&-s\end{pmatrix}}$
of left $\star$-multiplication by $\begin{pmatrix}1&-s(p)\end{pmatrix}$ is a contraction
from $\mathcal M$ into $\mathcal H(s)$, where $\mathcal H(s)$
denotes the reproducing kernel Hilbert space of quaternionic
valued
functions which are hyperholomorphic  in the ball $\mathbb B_1$ and with reproducing kernel
$$
K_s(p,q)=\sum_{t=0}^{\infty} p^t(1-s(p)\overline{s(q)})\bar q^t.
$$}\\

STEP 3: {\sl Assume that $s$ is a solution of Problem \ref{pb1} and that $P>0$. Then, $s$ is of the form \eqref{lft1}.}\\

STEP 4: {\sl Assume that $P>0$. Then any function of the form \eqref{lft1} satisfies
the interpolation condition \eqref{inter1} and if, in addition, \eqref{limits} is in force, then $s$ satisfies \eqref{inter2}}.\\

The proofs of Steps 1-4 are given in Section \ref{sec123}.
The degenerate case is considered in Section \ref{sec345}.

\section{Proofs of Steps 1-4}
\setcounter{equation}{0}
\label{sec123}

{\bf Proof of Step 1:} Assume a solution $s$ exists. Since $s$ is a Schur function the kernel $K_s(p,q)$ is positive definite and so for every $r\in(0,1)$ the $N\times N$ matrix $P(r)$ with $(u,v)$ entry equal to
\[
P_{uv}(r)=K_s(rp_u,rp_v)=\sum_{t=0}^\infty r^{2t}p_u^t(1-s(rp_u)\overline{s(rp_v)})p_v^t,\quad u,v=1,\ldots N
\]
is positive. Setting
\[
G=(1-s(rp_u)\overline{s(rp_v)}),\quad p=r^2p_u,\quad\text{and}\quad A=\overline{p_v}
\]
in formula \eqref{magic} we have
\[
P_{uv}(r)=\left((1-s(rp_u)\overline{s(rp_v)})-r^2\overline{p_u}(1-s(rp_u)\overline{s(rp_v)})\overline{p_v}\right)
(1-2r^2{\rm Re}(p_u)\overline{p_v}+r^4\overline{p_v}^2)^{-1}.
\]
Furthermore, we note that $P(r)$ is a solution of the matrix equation
\[
P(r)-r^2A^*P(r)A=C(r)^*\mathcal{J}C(r)
\]
where
\[
C(r)=\begin{pmatrix}1&\cdots &1\\
& &\\
\overline{s(rp_1)}&\cdots &\overline{s(rp_N)}\end{pmatrix},
\]
and $A$ is as in (\ref{Dn:ACJ}).
In fact, with the above notation, the $(u,v)$ element of the matrix $P(r)-r^2A^*P(r)A$ can be computed as follows:
\[
\begin{split}
&P_{uv}(r)-r^2p_uP_{uv}(r)\overline{p_v}\\
&=\left(
\left(G-r^2\overline{p_u}G\overline{p_v}\right) -r^2p_u\left(G-r^2\overline{p_u}G\overline{p_v}\right)\overline{p_v}
\right)
(1-2r^2{\rm Re}(p_u)\overline{p_v}+r^4\overline{p_v}^2)^{-1}\\
&=\left(
G-r^2\overline{p_u}G\overline{p_v} -r^2p_uG \overline{p_v}+r^4G\overline{p_v}^2
\right)
(1-2r^2{\rm Re}(p_u)\overline{p_v}+r^4\overline{p_v}^2)^{-1}\\
&=G\left(
1-2r^2{\rm Re}({p_u})\overline{p_v} +r^4\overline{p_v}^2
\right)
(1-2r^2{\rm Re}(p_u)\overline{p_v}+r^4\overline{p_v}^2)^{-1}=(1-s(rp_u)\overline{s(rp_v)})\\
\end{split}
\]
and so the $(u,v)$ element in the matrix $P(r)-r^2A^*P(r)A$ equals the $(u,v)$ element in $C(r)^* \mathcal{J}C(r)$ as stated.
We now let $r$ tend to $1$. Since $s$ is assumed to be a solution of Problem \ref{pb1}, we have
\[
\lim_{\substack{r\rightarrow 1\\
r\in(0,1)}}K_s(rp_u,rp_u)=\lim_{\substack{r\rightarrow 1\\
r\in(0,1)}}\frac{1-|s(rp_u)|^2}{1-r^2} \le \kappa_u,\quad
u=1,\ldots N
\]
and
\[
\lim_{\substack{r\rightarrow 1\\ r\in(0,1)}}C(r)=C,
\]
where $C$ is as in (\ref{Dn:ACJ}).
Furthermore we note that $1-2{\rm Re}(p_u)\overline{p_v}+\overline{p_v}^2\not=0$ since $1-2{\rm Re}(p_u)x+x^2$ is the so-called minimal (or companion) polynomial associated with the sphere $[p_u]$ which vanishes exactly at points on the sphere $[p_u]$ and $p_v\not\in [p_u]$. This fact can also be obtained directly using Lemma \ref{La:ph-hq}. Indeed, for indices $u\not= v$, we have
\begin{equation}
\label{aveiro:2014}
1-2{\rm Re} ( p_u)\overline{p_v}+\overline{p_v}^2=p_u(\overline{p_u}-\overline{p_v})-(\overline{p_u}-\overline{p_v})
\overline{p_v}\not=0,
\end{equation}
since $p_u$ and $p_v$ (and hence $p_u$ and $\overline{p_v}$) are assumed on different spheres for $u\not=v$. It follows that
$\lim_{\substack{r\rightarrow 1\\ r\in(0,1)}}
P_{uv}(r)$ exists and is in fact equal to $P_{uv}$ for $u\not =v$ by uniqueness of the solution of the equation
\begin{equation}
\label{secretequation}
x-p_ux\overline{p_v}=0.
\end{equation}
Hence $P\ge 0$ since $P(r)\ge 0$ for all $r\in(0,1)$.\\

{\bf Proof of Step 2:}  Let $s$ be a solution (if any) of Problem \ref{pb1}, let $u\in\left\{1,\ldots, N\right\}$,
and let $r\in(0,1)$. The functions
\[
g_{u,r}(p)=K_s(p, rp_u)=\sum_{t=0}^\infty
p^t(1-s(p)\overline{s(rp_u)})\overline{p_u}^t
\]
belong to $\mathcal H(s)$ and have uniformly bounded norms since
\[
\lim_{\substack{r\rightarrow 1\\
r\in(0,1)}}\|g_{u,r}(rp_u)\|^2_{\mathcal H(s)}=
\lim_{\substack{r\rightarrow 1\\
r\in(0,1)}}K_s(rp_u,rp_u)\le\kappa_u.
\]
Thus there is a sequence of numbers $r_0,r_1,\ldots \in(0,1)$ which tends to $1$ (without loss of generality we may assume that the sequence is the same for $p_1, \cdots, p_N$) and an element $g_u\in\mathcal H(s)$ such that the functions
$g_{u,r_n}$ tend weakly to $g_u$. In a reproducing kernel Hilbert space weak convergence implies pointwise convergence, and so
\[
\begin{split}
g_u(p)&=\lim_{n\longrightarrow\infty}g_{u,r_n}(p)\\
&=\lim_{n\longrightarrow\infty}
\sum_{t=0}^\infty r_n^tp^t(1-s(p)\overline{s(r_np_u)})\overline{p_u}^t\\
&=\sum_{t=0}^\infty p^t(1-s(p)\overline{s_u})\overline{p_u}^t\\
&=\begin{pmatrix}1&-s(p)\end{pmatrix}\star f_u(p),\quad\forall p\in\mathbb B_1,
\end{split}
\]
where
\begin{equation}
\label{fu}
f_u(p)=\sum_{t=0}^\infty p^t\begin{pmatrix}1\\ \overline{s_u}\end{pmatrix}\overline{p_u}^t
\end{equation}
denotes the $u$-th column of the matrix-function $F(p)$ and where
the interchange of summation and limit is justified since $|p|<1$. Hence $M_{\begin{pmatrix}1&-s\end{pmatrix}}$ sends
$\mathcal M$ into $\mathcal H(s)$.
Note that for $Y = (y_{u,v})_{u,v=1}^N$ and $Z = (z_{u,v} )_{u,v=1}^N$ we define $Y \star Z$ to be the $N \times N$ matrix whose $(u,v)$ entry is given by $\sum_{t=1}^N y_{u,t} \star z_{t,v}$.
To show that this operator is a contraction we first compute the inner product $\langle g_v,g_u
\rangle_{\mathcal H(s)}$ for $u\not= v$. By the definition of the weak limit and of the reproducing kernel, we can write
\[
\begin{split}
\langle g_v,g_u
\rangle_{\mathcal H(s)}&=\lim_{n\longrightarrow\infty}\langle g_v,g_{u,r_n}
\rangle_{\mathcal H(s)}\\
&=\lim_{n\longrightarrow\infty}g_v(r_np_u)\\
&=\lim_{n\longrightarrow\infty}\sum_{t=0}^\infty r_n^tp_u^t(1-s(r_np_u)\overline{s_v})\overline{p_v}^t\\
&=\lim_{n\longrightarrow\infty}\left((1-s(r_np_u)\overline{s_v})-r_n\overline{p_u}(1-s(r_np_u)\overline{s_v})\overline{p_v}\right)
(1-2r_n{\rm Re}( p_u)\overline{p_v}+r_n^2\overline{p_v}^2)^{-1}\\
&=\left((1-s_u\overline{s_v})-\overline{p_u}(1-s_u\overline{s_v})\overline{p_v}\right)
(1-2{\rm Re}( p_u)\overline{p_v}+\overline{p_v}^2)^{-1},
\end{split}
\]
where we have used formula \eqref{magic} and, as in the proof of Step 1
(see \eqref{aveiro:2014}), the fact that $[p_u]\cap[p_v]=\emptyset$ (recall that we assume here $u\not=v$). We claim that
\begin{equation}
\label{parisbastille}
P_{uv}=\left((1-s_u\overline{s_v})-\overline{p_u}(1-s_u\overline{s_v})\overline{p_v}\right)
(1-2{\rm Re} (p_u)\overline{p_v}+\overline{p_v}^2)^{-1}.
\end{equation}
The proof is similar to the argument in the proof of step 1, and is as follows. Set $h_n=\langle g_v,g_{u,r_n}
\rangle_{\mathcal H(s)}$. Then
\[
h_n-r_np_uh_n\overline{p_v}=1-s(r_np_u)\overline{s_v}.
\]
Letting $n\rightarrow\infty$ we see that $h=\lim_{n\rightarrow\infty}h_n$ satisfies equation \eqref{puv}.
By the uniqueness of the solution of this equation we have $h=P_{uv}$. Furthemore, by the property of the weak limit versus the norm,
\begin{equation}
\label{parisstmichel}
\|g_u\|^2_{\mathcal H(s)}\le\lim_{n\rightarrow\infty}\|g_{u,r_n}\|^2_{\mathcal H(s)}\le \kappa_u.
\end{equation}
We can now show that $\|M_{\begin{pmatrix}1&-s\end{pmatrix}}\|\le 1$.
Let $c\in \mathbb H^N$. Then,
\[
\left(M_{\begin{pmatrix}1&-s\end{pmatrix}}Fc\right)(p)=\sum_{u=1}^N g_u(p)c_u
\]
and we have
\[
\begin{split}
\|(M_{\begin{pmatrix}1&-s\end{pmatrix}}Fc\|^2_{\mathcal H(s)}&=\sum_{u,v=1}^N
\overline{c_u}\left(\langle g_v,g_u\rangle_{\mathcal H(s)}\right)c_v\\
&=\sum_{u=1}^N
|c_u|^2\|g_u\|^2_{\mathcal H(s)}+\sum_{\substack{u,v=1\\ u\not=v}}^N
\overline{c_u}\left(\langle g_v,g_u\rangle_{\mathcal H(s)}\right)c_v\\
&=\sum_{u=1}^N
|c_u|^2\|g_u\|^2_{\mathcal H(s)}+\sum_{\substack{u,v=1\\ u\not=v}}^N
\overline{c_u}P_{uv}c_v\\
&\le \sum_{u=1}^N|c_u|^2\kappa_u+\sum_{\substack{u,v=1\\ u\not=v}}^N
\overline{c_u}P_{uv}c_v\\
&=c^*Pc\\
&=\|Fc\|^2_{\mathcal M},
\end{split}
\]
where we have used \eqref{parisbastille} and \eqref{parisstmichel}. Thus the $\star$-multiplication by $(1\ -s(p))$ is a contraction from $\mathcal M$ into $\mathcal{H}(s)$.\\

{\bf Proof of Step 3:} Let $\Theta$ be defined by \eqref{parislelouvre}, and
\begin{equation}
\label{musee-d-orsay}
K_\Theta(p,q)=\sum_{t=0}^\infty
p^t\left(\mathcal{J}-\Theta(p)\mathcal{J}\Theta(q)^*\right)\overline{q}^t.
\end{equation}
The formula
\begin{equation}
\label{fpf}
F(p)P^{-1}F(q)^*=K_\Theta(p,q)
\end{equation}
is proved as in the complex case when $p$ and $q$ are real, and is then extended to $p,q\in\mathbb B_1$ by a slice hyperholomorphic
extension. Using \eqref{grenelle} we have
\[
\left(M_{\begin{pmatrix}1&-s\end{pmatrix}}^*K_s(\cdot, q)\right)(p)=\sum_{t=0}^\infty
p^n\left(\begin{pmatrix}1\\ -\overline{s(q)}\end{pmatrix}-\Theta(p)\mathcal{J}\Theta(q)^*\star_r\begin{pmatrix}1\\
-\overline{s(q)}\end{pmatrix}\right)\overline{q}^t,
\]
and so
\[
\begin{split}
(M_{\begin{pmatrix}1&-s\end{pmatrix}}M_{\begin{pmatrix}1&-s\end{pmatrix}}^*K_s(\cdot, q))(p)&\\
&\hspace{-2cm}=K_s(p,q)-
\sum_{t=0}^\infty p^t\left(\begin{pmatrix}1&-s(p)\end{pmatrix}\star\Theta(p)
\mathcal{J}\Theta(q)^*\star_r\begin{pmatrix}1\\
-\overline{s(q)}\end{pmatrix}\right)\overline{q}^t\\
&\hspace{-2cm}\le K_s(p,q),
\end{split}
\]
and therefore the kernel
\[
\sum_{t=0}^\infty p^t\left(\begin{pmatrix}1&-s(p)\end{pmatrix}\star\Theta(p)
\mathcal{J}\Theta(q)^*\star_r\begin{pmatrix}1\\
-\overline{s(q)}\end{pmatrix}\right)\overline{q}^t
\sum_{t=0}^\infty p^t\left(A(p)\overline{A(q)}-B(p)\overline{B(q)}\right)\overline{q}^t
\]
is positive definite in $\mathbb B_1$, where
\[
A(p)=(a-s\star c)(p)\quad\text{and}\quad  B(p)=(b-s\star d)(p).
\]
The point $p=1$ is not an interpolation node, and so $\Theta$ is well defined at $p=1$. From \eqref{parislelouvre} we have
\begin{equation}
\label{operabastille}
\Theta(1) = \mathcal{I}_2
\end{equation}
and so $(a^{-1}c)(1) = 0$. Since $s$ is bounded by $1$ in modulus in $\mathbb B_1$ it follows that $(a-s\star c)(p)\not\equiv 0$. Thus
$e=-(a-s\star c)^{-\star}\star(b-s\star d)$ is defined in $\mathbb B_1$, with the possible exception of spheres of poles. Since
\[
\sum_{t=0}^\infty p^t\left(A(p)\overline{A(q)}-B(p)\overline{B(q)}\right)\overline{q}^t=A(p)\star\left\{
\sum_{t=0}^\infty p^t(
1-e(p)\overline{e(q)})\overline{q}^t\right\}\star_r
\overline{A(q)},
\]
we have from \cite[Proposition 5.3]{acs3} that the kernel
\[
K_e(p,q)=\sum_{t=0}^\infty p^t(
1-e(p)\overline{e(q)})\overline{q}^t
\]
is positive definite in its domain of definition, and thus $e$ extends to a Schur function (see \cite{acs1} for the latter assertion). From
\[
e=-(a-s\star c)^{-\star}\star(b-s\star d)
\]
we get $s\star(c\star e +d)=a\star e+b$. To conclude we remark that \eqref{operabastille} implies that
\[
(d^{-1}c)(1)=0.
\]
Thus, as just above $c\star e+d\not\equiv 0$ and we get  that $s$ is of the form \eqref{lft1}.\\

{\bf Proof of Step 4:} Assume that $s$ is of the form \eqref{lft1}. Then the formula
\[
K_s(p,q)=\begin{pmatrix}1&-s(p)\end{pmatrix}\star K_\Theta(p,q)\star_r \begin{pmatrix}1\\-\overline{s(q)}\end{pmatrix}+
(a-s\star c)(p)\star K_e(p,q)\star_r\overline{(a-s\star c)(q)}
\]
implies that $M_{\begin{pmatrix}1&-s\end{pmatrix}}$ is a contraction from $\mathcal H(\Theta)$ into $\mathcal H(s)$. In particular
\begin{equation}
\label{gu}
g_u(p)=\begin{pmatrix}1&-s(p)\end{pmatrix}\star f_u(p)=\sum_{t=0}^\infty p^t(1-s(p)\overline{s_u})\overline{p_u}^t\in\mathcal H(s)
\end{equation}
and
\[
\|g_u\|^2_{\mathcal H(s)}\le \kappa_u.
\]
We want to infer from these facts that $s$ satisfies the
interpolation conditions \eqref{inter1}. We have
\begin{equation}
\label{richelieu-drouot}
\begin{split}
|g_u(rp_u)|^2&=|\langle g_u(\cdot), K_s(\cdot, rp_u)\rangle_{\mathcal H(s)}|^2\\
&\le \left(\|g_u\|^2_{\mathcal H(s)}\right)\cdot K_s(rp_u,rp_u)\\
&\le \kappa_u\cdot\frac{1-|s(rp_u)|^2}{1-r^2}\\
&\le \frac{2\kappa_u}{1-r}.
\end{split}
\end{equation}
In view of \eqref{magic}, we get
\begin{equation}
\label{republique}
\begin{split}
g_u(rp_u)&=\sum_{t=0}^\infty r^tp_u^t(1-s(rp_u)\overline{s_u})\overline{p_u}^t\\
&=\left((1-s(rp_u)\overline{s_u})-r\overline{p_u}(1-s(rp_u)\overline{s_u})\overline{p_u}\right)(1-2r{\rm Re}( p_u)
\overline{p_u}+r^2\overline{p_u}^2)^{-1}\\
&=\left((1-s(rp_u)\overline{s_u})-r\overline{p_u}(1-s(rp_u)\overline{s_u})\overline{p_u}\right)((1-r)(1-r\overline{p_u}^2))^{-1},
\end{split}
\end{equation}
and so we have
\[
\frac{|(1-s(rp_u)\overline{s_u})-r\overline{p_u}(1-s(rp_u)\overline{s_u})\overline{p_u}|}{|1-r\overline{p_u}^2|}\le \sqrt{2\kappa_u}\cdot
\sqrt{1-r}.
\]
Let $\sigma_u$ be a limit, via a subsequence, of $s(rp_u)$ as $r\rightarrow 1$, and set $X_u=1-\sigma_u\overline{s_u}$.
The above inequality implies that $X_u=\overline{p_u}X_u\overline{p_u}$,
and so
\begin{equation}
\label{richardlenoir}
X_up_u=\overline{p_u}X_u.
\end{equation}
The conjugate of \eqref{richardlenoir} is
\begin{equation}
\label{richardlenoir1}
\overline{X_u}p_u=\overline{p_u}\overline{X_u}.
\end{equation}
Adding \eqref{richardlenoir} and \eqref{richardlenoir1} we obtain
\[
{\rm Re}( X_u)p_u=\overline{p_u}{\rm Re}( X_u).
\]
Since $p_u$ is not real we get that ${\rm Re}( X_u)=0$. Let
$X_u=\alpha i+\beta j+\gamma k$, where $\alpha, \beta,\gamma\in\mathbb R$. From
$\sigma_u\overline{s_u}=1-X_u$ we have
\[
|\sigma_u\overline{s_u}|^2=1+\alpha^2+\beta^2+\gamma^2.
\]
Since $\sigma_u\in\mathbb B_1$ we have $|\sigma_u\overline{s_u}|\le 1$ and so $\alpha=\beta=\gamma=0$. Thus, $X_u=0$ and $\sigma_u\overline{s_u}=1$. Hence
$\sigma_u=s_u$ and the limit $\lim_{\substack{r\rightarrow 1\\ r\in(0,1)}} s(rp_u)$ exists and is equal to $s_u$, and hence \eqref{inter1} is satisfied.\\

To prove that \eqref{inter2} is met we proceed as follows. From  \eqref{richelieu-drouot} we have in particular
\[
|g_u(rp_u)|^2
\le \kappa_u\cdot\frac{1-|s(rp_u)|^2}{1-r^2},
\]
and using \eqref{republique} we obtain:
\begin{equation}
\label{jussieu}
\frac{|X(r)-r\overline{p_u}X(r)\overline{p_u}|^2}{(1-r)^2|1-r\overline{p_u}^2|^2}\le
\kappa_u\cdot\frac{1-|s(rp_u)|^2}{1-r^2},
\end{equation}
where we have set $X(r)=1-s(rp_u)\overline{s_u}$. Assume now that
\eqref{limits} is in force and let
\begin{equation}
\lim_{\substack{r\rightarrow 1\\
r\in(0,1)}}\frac{1-s(rp_u)\overline{s_u}}{1-r}=\beta_u\in\mathbb
R.
\end{equation}
Then \eqref{jussieu} together with \eqref{richelieu-drouot1}
imply that
\[
\beta_u^2\le \beta_u\kappa_u,
\]
from which we get that $\beta_u\ge0$ and
\[
\lim_{\substack{r\rightarrow 1\\
r\in(0,1)}} \frac{1-s(rp_u)\overline{s_u}}{1-r}\le \kappa_u.
\]
\section{The degenerate case}
\setcounter{equation}{0}
\label{sec345}

We now consider the case where $P$ is singular. We need first a
definition. A finite Blaschke product is a finite $\star$-product of terms of the form which are given by
\begin{equation}
\label{eqBlaschke} b_a(p)=(1-p\bar
a)^{-\star}\star(a-p)\frac{\bar a}{|a|},
\end{equation}
where $a\in\mathbb{H}$, $|a|<1$ (see \cite{MR3127378}).\\

The purpose of this section is to prove the following theorem.
First a remark. We denote by $r$ the rank of $P$ and assume that
the main $r\times r$ minor of $P$ is invertible. This can be done
by rearranging the interpolation points.

\begin{Tm}
\label{tm2} Assume that $P$ is singular. Then Problem \ref{pb1}
has at most one solution, and the latter is then a finite
Blaschke product. It has a unique solution satisfying
\eqref{bastille} for $u=1,\ldots, r$.
\end{Tm}

We begin with some preliminary results and definitions.

\begin{Dn}
Let $f$ be a slice hyperholomorphic in a neighborhood $\Omega$ of $p=1$, and let $f(p)=\sum_{t=0}^\infty (p-1)^tf_t$
be its power series expansion at $p=1$. We define
\begin{equation}
R_1f(p)=\sum_{t=1}^\infty (p-1)^t f_t.
\end{equation}
\end{Dn}

Denoting by ${\rm ext}$ the slice hyperholomorphic extension we have
\begin{equation}
R_1f(p)={\rm ext}\left(R_1f|_{p=x}\right).
\label{lisbonne}
\end{equation}
\begin{La}\label{lemma:R1}
Let $f(p)=F(p)\xi$ where $F(p)=C\star (\mathcal{I}_N-pA)^{-\star}$, then
\begin{equation}
\label{place-des-vosges}
R_1f(p)=F(p)A(\mathcal{I}_N-A)^{-1}\xi.
\end{equation}
\end{La}
\begin{proof}
First of all, recall that 
$$
F(p)=C\star (\mathcal{I}_N-pA)^{-\star}=(C-\bar p CA)(I_n-2{\rm Re}(p) A +|p|^2 A^2)^{-1}
$$
so 
$$
F(1)=(C-CA)(\mathcal{I}_N-2 A + A^2)^{-1}=C(\mathcal{I}_N-A)^{-1}.
$$
Let us compute
\[
\begin{split}
R_1 f(p)&=(p-1)^{-1}(f(p)-f(1))=(p-1)^{-1}(C\star (\mathcal{I}_N-pA)^{-\star}\xi- C (\mathcal{I}_N-A)^{-1}\xi)\\
&= C\star (p-1)^{-1}( (\mathcal{I}_N-pA)^{-\star}- (\mathcal{I}_N-A)^{-1})\xi\\
&= C\star (p-1)^{-1} \star (\mathcal{I}_N-pA)^{-\star}\star ((\mathcal{I}_N-A) - (\mathcal{I}_N-pA))(\mathcal{I}_N-A)^{-1}\xi\\
&= C\star (p-1)^{-1} \star (\mathcal{I}_N-pA)^{-\star}\star (p-1) A(\mathcal{I}_N-A)^{-1}\xi\\
&= C\star  (\mathcal{I}_N-pA)^{-\star} A(\mathcal{I}_N-A)^{-1}\xi \\
&= F(p)A(\mathcal{I}_N-A)^{-1}\xi.
\end{split}
\]
\end{proof}
\begin{La}
Let $f,g\in\mathcal M$. Then
\begin{equation}
\label{dbspecial-case} [f,g]+[R_1f,g]+[f,R_1g]=g(1)^*\mathcal{J} f(1).
\end{equation}
\end{La}
\begin{proof} Let $f(p)=F(p)\xi$ and $g(p)=F(p)\eta$ with $\xi,\eta\in\mathbb H^N$. We have
\[
f(1)=C(\mathcal{I}_N-A)^{-1}\xi \quad\text{and }\quad g(1)=C(\mathcal{I}_N-A)^{-1}\eta.
\]
These equations together with \eqref{place-des-vosges} show that \eqref{dbspecial-case}
is equivalent to
\[
P+P(\mathcal{I}_N-A)^{-1}A+A^*(\mathcal{I}_N-A)^{-*}P=(\mathcal{I}_N-A)^{-*}C^*\mathcal{J}C(\mathcal{I}_N-A).
\]
Multiplying this equation by $\mathcal{I}_N-A^*$ on the left and by $\mathcal{I}_N-A$ on the right we get the equivalent equation
\eqref{steineq}.
\end{proof}

\begin{Rk}
{\rm
Equation \eqref{dbspecial-case} corresponds to a special case of a structural identity which characterizes $\mathcal H(\Theta)$ spaces in the complex setting. A corresponding identity in the half place case was first introduced by
de Branges, see \cite{dbhsaf1}, and improved by Rovnyak \cite{HM}. Ball introduced the corresponding identity in the setting of the
open unit disk and proved the corresponding structure theorem. See \cite{ball-contrac}.
See e.g.  \cite[p. 17]{ad-jfa} for further discussions on this topic.
}
\end{Rk}

\begin{Pn}
Let $a$ and $b$ be slice hyperholomorphic functions defined in an axially symmetric s-domain containing $p=1$. Then,
\begin{equation}
\label{opera-garnier}
R_1(a\star b)(p)=\left(R_1a(p)\right)b(1)+(a\star R_1b)(p) .
\end{equation}
\end{Pn}

\begin{proof}
By the Identity Principle, see \cite[Theorem 4.2.4]{MR2752913} the equality holds if and only if it holds for the restrictions to a complex plane $\mathbb C_I$ i.e., using the notations in Section 2, if and only if
\begin{equation}\label{R1equality}
(R_1(a\star b))_{I}(z)=
\left(R_1a(z)\right)_{I}b(1)+(a\star R_1b)_{I}(z), \quad z\in\mathbb C_I.
\end{equation}
Let $J\in\mathbb S$ be such that $J$ is orthogonal to $I$ and assume that
$$
a_{I}(z)=F(z)+G(z)J,\qquad b_{I}(z)=H(z)+L(z)J .
$$
Let us compute the left-hand side of (\ref{R1equality}), using the fact that $(R_1(a\star b))_{I}(z)= R_1((a\star b)_{I})$ and formula (\ref{starproduct}):
\[
\begin{split}
R_1((a\star b)_{I})&=R_1\left(F(z)H(z)-G(z)\overline{L(\bar z)}+(G(z)\overline{H(\bar z)}+F(z)L(z))J\right) \\
&=(z-1)^{-1}\left(F(z)H(z)-G(z)\overline{L(\bar z)}+(G(z)\overline{H(\bar z)}+F(z)L(z))J\right.\\
&\left.-F(1)H(1)+G(1)\overline{L(1)}-(G(1)\overline{H(1)}+F(1)L(1))J)\right).
\end{split}
\]
At the right hand side of (\ref{R1equality}) we have $\left(R_1a(z)\right)_{I}b(1)=\left(R_1a_{I}(z)\right)b(1)$ which can be written as
\[
\begin{split}
\left(R_1a_{I}(z)\right)&b(1)=\left((z-1)^{-1}(F(z)+G(z)J-F(1)-G(1)J)\right)(H(1)+L(1)J)\\
&=(z-1)^{-1}\left(F(z)H(1)+F(z)L(1)J+G(z)\overline{H(1)}J-G(z)\overline{L(1)}-F(1)H(1)\right.\\
&\left.-F(1)L(1)J-G(1)\overline{H(1)}J+G(1)\overline{L(1)}\right),
\end{split}
\]
moreover
\[
\begin{split}(a\star R_1b)_{I}(z)&=(F(z)+G(z)J)\star\left((z-1)^{-1}(H(z)+L(z)J-H(1)-L(1)J)\right)\\
&= (z-1)^{-1} (F(z)+G(z)J) \star (H(z)+L(z)J-H(1)-L(1)J)\\
&=(z-1)^{-1}(F(z)H(z)-G(z)\overline{L(\bar z)}+(G(z)\overline{H(\bar z)}+F(z)L(z))J )\\
& -F(z)H(1)+G(z)\overline{L(1)}-(G(z)\overline{H(1)}+F(z)L(1))J
\end{split}
\]
from which the equality follows.
\end{proof}

We will also need the following result, well known in the complex case. We refer to \cite{aron,schwartz} for more information and to
\cite{fw} for connections with operator ranges.

\begin{Tm}
\label{new}
Let $K_1(p,q)$ and $K_2(p,q)$ be two $\mathbb H$-valued functions positive definite in a set $\Omega$ and assume that
the corresponding reproducing kernel Hilbert spaces have a zero intersection. Then the sum
\[
\mathcal H(K_1+K_2)=\mathcal H(K_1)+\mathcal H(K_2)
\]
is orthogonal.
\end{Tm}

\begin{proof}  Let $K=K_1+K_2$. The linear relation in $\mathcal H(K)\times (\mathcal H(K_1)\times \mathcal H(K_2))$ spanned
by the pairs
\[
(K(p,q), (K_1(p,q),K_2(p,q))),\quad q\in\Omega,
\]
is densely defined and isometric. It therefore extends to the graph of an everywhere defined isometry, which we will call $T$. See
\cite[Theorem 7.2]{MR3127378}. From
\[
\begin{split}
(T^*(f_1,f_2))(q)&=\langle T^*(f_1,f_2),K(p,q)\rangle_{\mathcal H(K)}\\
&=\langle (f_1,f_2), TK(p,q)\rangle_{\mathcal H(K_1)\times \mathcal H(K_2)}\\
&=\langle f_1, K_1(p,q)\rangle_{\mathcal H(K_1)}+\langle f_2, K_2(p,q)\rangle_{\mathcal H(K_2)}\\
&=f_1(q )+f_2(q),\quad q\in\Omega,
\end{split}
\]
we see that $\ker T^*=\left\{0\right\}$ since $\mathcal H(K_1)\cap \mathcal H(K_2)=\left\{0\right\}$. Thus $T$ is unitary and the result
follows then easily.
\end{proof}

\begin{proof}[Proof of Theorem \ref{tm2}]
We proceed in a number of steps. Recall that $r={\rm rank}\, P$.\\

STEP 1: {\sl Assume $r=0$. Then, $s_1=\cdots=s_N$ and Problem \ref{pb1} is solvable with the unique solution
the constant unitary function $s(p)\equiv s_1$ .}\\

The matrix $P=0$, and equation \eqref{steineq} imply that
$C^*\mathcal{J}C=0$, and so $1-s_u\overline{s_v}=0$ for $u\not= v\in\left\{
1,\ldots, N\right\}$. Thus $s_1=\cdots=s_N$ and the function
$s(p)\equiv s_1$ is clearly a solution. Assume that $s$ is a
(possibly different) solution of Problem \ref{pb1}. The
map $M_{\begin{pmatrix}1&-s\end{pmatrix}}$ of slice
multiplication by $\begin{pmatrix}1&-s(p)\end{pmatrix}$ is a
contraction from $\mathcal M$ into $\mathcal H(s)$ (see the
second step in the proof of Theorem \ref{tm1}). Thus
\[
\begin{pmatrix}1&-s(p)\end{pmatrix}\star f_u(p)\equiv 0,\quad u=1,\ldots, N,
\]
that is $g_u(p)\equiv 0$, where $f_u$ and $g_u$ have been defined in \eqref{fu} and \eqref{gu} respectively.
From \eqref{magic} we have (for $|p|<1$)
\[
g_u(p)=\left((1-s(p)\overline{s_u})-\overline{p}(1-s(p)\overline{s_u})\overline{p_u}\right)(1-2{\rm Re}(p)\overline{p_u}+
|p|^2p_u^2)^{-1}
\]
since
\[
1-2{\rm Re}(p)\overline{p_u}+|p|^2p_u^2\not=0
\]
for $|p|<1$. Hence
\[
(1-s(p)\overline{s_u})=\overline{p}(1-s(p)\overline{s_u})\overline{p_u},\quad\forall p\in\mathbb H_1.
\]
Taking absolute values of both sides of this equality we get $1-s(p)\overline{s_u}\equiv 0$, and so $s(p)\equiv s_u$.
This ends the proof of Step 1.\\

In the rest of the proof we assume $r > 0$. By reindexing the
interpolating nodes we can assume that the principal minor of
order $r$ is invertible. Thus the corresponding space is a
$\mathcal H(\Theta_r)$ space, and we can write
\[
\mathcal M=\mathcal H(\Theta_r)\oplus\Theta_r\star \mathcal N.
\]

STEP 2: {\sl The elements of $\mathcal N$ are slice hyperholomorphic in a
neighborhood of $p=1$ and
$R_1\mathcal N\subset\mathcal N$.}\\

We follow the argument in Step 1 in the proof of Theorem 3.1 in
\cite{ad-laa-herm} (see p. 153). From \eqref{opera-garnier} we
have
\begin{equation}
(R_1(\Theta_r\star n))(p)=(R_1\Theta_r)(p)n(1)+(\Theta_r\star
R_1n)(p). \label{newlabel}
\end{equation}
To prove that $R_1n\in\mathcal N$ we show that
\begin{equation}
\label{au-clair-de-la-lune}
[(R_1(\Theta_r\star
n))(p)-(R_1\Theta_r)(p)n(1),g]_{\mathcal M}=0,\quad \forall
g\in\mathcal H(\Theta_r).
\end{equation}
Using \eqref{dbspecial-case} we have
\[
\begin{split}
[(R_1(\Theta_r\star n))(p),g]_{\mathcal
M}&=g(1)^*\mathcal{J}(R_1(\Theta_r\star n))(1)-[\Theta_r\star
n,g]_{\mathcal M}-[\Theta_r\star n,
R_1g]_{\mathcal M}\\
&=g(1)^*\mathcal{J}(R_1(\Theta_r\star n))(1)
\end{split}
\]
since
\[
[\Theta_r\star n,g]_{\mathcal M}=0\quad\text{and}\quad
[\Theta_r\star n, R_1g]_{\mathcal M}=0,
\]
where the second equality follows from $R_1g\in\mathcal M$.
Moreover, for real $p=x$ we have the equality of real analytic functions
\[
(R_1\Theta_r)(x)=-K_{\Theta_r}(x,1)\mathcal{J}\Theta_r(1)^*,
\]
and so, by slice hyperholomorphic extension, see  \cite[Remark 2.18]{2013arXiv1310.1035A}, in a suitable neighborhood of $p=1$ we have
\[
(R_1\Theta_r)(p)=-K_{\Theta_r}(p,1)\mathcal{J}\Theta_r(1)^*.
\]
Note that $\Theta_r(1)$ is the identity.
Thus
\[
\begin{split}
[(R_1\Theta_r)(p)n(1),g]_{\mathcal
M}&=-[K_{\Theta_r}(p,1)\mathcal{J}\Theta_r(1)^*n(1),g]_{\mathcal M}\\
&=-(n(1)^*\Theta_r(1)^*g(1)^*)\\
&=-g(1)^*\Theta_r(1)\mathcal{J}n(1),
\end{split}
\]
and so \eqref{au-clair-de-la-lune} is in force. This ends the
proof of the second step.\smallskip

Endow now $\mathcal N$ with the Hermitian form
\[
[n_1, n_2]_{\mathcal N}= [\Theta_r\star n_1,\Theta_r\star
n_2]_{\mathcal M}.
\]

STEP 3: {\sl There exist matrices $(G,T)\in\mathbb H^{2\times
(N-r)}\times\mathbb H^{(N-r)\times (N-r)}$ such that $\mathcal N$
is spanned by the columns of the function $F_{\mathcal
N}(p)=G\star(\mathcal{I}_{N-r}-pT)^{-\star}$ and moreover for
$\xi\in\mathbb H^{N-r}$.
\[
F_{\mathcal N}(p)\xi\equiv 0\quad\Longrightarrow\quad \xi=0.
\]
}

Indeed, we first note that the elements of $\mathcal N$ are well
defined at $p=1$ since $\Theta$ is invertible at $p=1$ (see also
the formulas in \cite[Theorem 3.3 (2)]{ad-laa-herm} ).
Let $F_{\mathcal N}(p)$ be built from the columns of a basis of
$\mathcal N$ and note that there exists $B\in\mathbb H^{(N-r)\times (N-r)}$
such that
\[
R_1F_{\mathcal N}=F_{\mathcal N}B.
\]
Restricting to $p =x$, where $x$ is real, we have
\[
\frac{F(x)-F(1)}{x-1}=F(x)B,
\]
and so
\begin{equation}
\label{le-louvre}
F(x)(\mathcal{I}_{N-r}+B-xB)=F(1).
\end{equation}
We claim that $\mathcal{I}_{N-r}+B$ is invertible. Let $\xi\in\mathbb H^{N-r}$ be
such that $B\xi=-\xi$. Then, \eqref{le-louvre} implies that
\[
xF(x)\xi=F(1)\xi,\quad x\in(-1,1).
\]
Thus $F(1)\xi=0$ (by setting $x=0$) and so $F(x)\xi=0$ and so
$\xi=0$. Hence
\[
F(x)=F(1)(\mathcal{I}_{N-r}+B)^{-1}(\mathcal{I}_{N-r}-xB(\mathcal{I}_{N-r}+B)^{-1})^{-1},
\]
and the result follows.\\

 The following step is \cite[Step 2 of proof of Theorem 3.1, p.
154]{ad-laa-herm}. The proof uses \eqref{au-clair-de-la-lune} and
is similar to the above arguments. \\

STEP 3: {\sl The space $\mathcal N$ is neutral and $G^*\mathcal{J}G=0$.} \\

%
%
$\mathcal N$ is neutral by construction since $r= {\rm rank} \, P$. We first show that the inner product in $\mathcal N$
satisfies \eqref{dbspecial-case}. We may proceed as in \cite[p.
154]{ad-laa-herm} and using \eqref{dbspecial-case} in $\mathcal
M$ we have for $n_1,n_2\in\mathcal M$: \[
\begin{split}
[R_1n_1,n_2]_{\mathcal N}&=[\Theta\star R_1n_1,\Theta\star
n_2]_{\mathcal M}\\
&=[R_1(\Theta\star n_1),\Theta\star n_2]_{\mathcal M}-
[(R_1\Theta)(n_1(1)),\Theta\star n_2]_{\mathcal M}\quad
\text{(where we used \eqref{newlabel})}\\
&=[R_1(\Theta\star n_1),\Theta\star n_2]_{\mathcal M}
\end{split}
\]
since $(R_1\Theta)(n_1(1))\in\mathcal H(\Theta)$, and so
$[(R_1\Theta)(n_1(1)),\Theta\star n_2]_{\mathcal M}=0$.\smallskip

Similarly,
\[
\begin{split}
[n_1,R_1n_2]_{\mathcal N}&= [\Theta\star n_1,\Theta\star
R_1n_2]_{\mathcal M}\\
&=[\Theta\star n_1,(R_1\Theta)(n_2(1))]_{\mathcal M}-
[\Theta\star n_1,(R_1\Theta)(n_2(1))]_{\mathcal M}\\
&=[\Theta\star n_1,(R_1\Theta)(n_2(1))]_{\mathcal M}.
\end{split}
\]
Thus, with $m_1=\Theta\star n_1$ and $m_2=\Theta\star n_2$,

\[
\begin{split}
[n_1,n_2]_{\mathcal N}+[R_1n_1,n_2]_{\mathcal N}+
[n_1,R_1n_2]_{\mathcal N}&=[m_1,m_2]_{\mathcal
M}+[R_1m_1,m_2]_{\mathcal M}+ [m_1,R_1m_2]_{\mathcal M}\\
&=m_2(1)^*\mathcal{J}m_1(1)\\
&=n_2(1)\mathcal{J}n_1(1)
\end{split}
\]
since $m_v(1)=(\Theta\star n_v)(1)=\Theta(1)n_v(1)$ for $v=1,2$ and
$\Theta(1)^*\mathcal{J}\Theta(1)=\mathcal{J}$.\\

Proceeding as in Step 1 it follows that
\[
P_{\mathcal N}-T^*P_{\mathcal N}T=G^*\mathcal{J}G,
\]
and so $G^*\mathcal{J}G=0$.\\

STEP 4: {\sl Problem \ref{pb1} has at most one solution.} \\

Let
\[
\Theta_r(p)=\begin{pmatrix}a_r(p)&b_r(p)\\c_r(p)&d_r(p)\end{pmatrix}.
\]
From the study of the nondegenerate case, we know that, under the assumptions that ensure the existence of a solution, any
solution is of the form
\begin{equation}
\label{ttheta} s(p)=(a_r(p)\star e(p)+b_r(p))\star(c_r(p)\star
e(p)+d_r(p))^{-\star},
\end{equation}
for some Schur function $e$. Furthermore as in step 1, for every
$n\in\mathcal N$ we have
\[
\begin{pmatrix}1&-s\end{pmatrix}\star\Theta_r\star n\equiv 0.
\]
Thus
\[
(a-sc)\star \begin{pmatrix}1&-e\end{pmatrix}\star n\equiv 0,
\]
and so
\[
\begin{pmatrix}1&-e\end{pmatrix}\star n\equiv 0.
\]
Since $G^*\mathcal{J}G=0$ we conclude in the way as in step 1. Indeed, let
\[
G=\begin{pmatrix}h_1 & \ldots &h_{N-r}\\k_1 & \ldots
&k_{N-r}\end{pmatrix}.
\]
At least one of the $h_u$ or $k_u$ is different from $0$ and
$G^*\mathcal{J}G=0$ implies that
\[
\overline{h_u}h_v=\overline{k_u}k_v,\quad\forall u,v=1,\ldots,
N-r,
\]
and so $e$ is a unitary constant.\\

We now show that the solution, when it exists, is a finite Blaschke product.\\

STEP 5: {\sl Let $s$ be given by \eqref{ttheta}. Then the
associated space $\mathcal H(s)$ is finite dimensional.}\smallskip

This follows from
\[
K_s(p,q)=\begin{pmatrix}1 &-s(p)\end{pmatrix}\star
K_{\Theta_r}(p,q)\star_r
\begin{pmatrix}1 \\ \overline{s(q)}\end{pmatrix}+
\underbrace{\begin{pmatrix}1 &-s(p)\end{pmatrix}\star
\Theta_r(p)\mathcal J\Theta_r(q)^* \star_r
\begin{pmatrix}1 \\ \overline{s(q)}\end{pmatrix}}_{\text{is equal to $0$ since
$|e|=1$}},
\]
where $K_{\Theta_r}$ is defined as in \eqref{musee-d-orsay} (with
$\Theta_r$ in place of $\Theta$).\\

STEP 6: {\sl The space $\mathcal H(s)$ contains an element of the
form
\begin{equation}
\label{sevres-babylone} f(p)=x\star(1-p\overline{a})^{-\star},
\end{equation}
where $x\in\mathbb H$ and $a\in\mathbb B_1$.}\smallskip

We first recall that (see \cite[Theorem 7.1]{acs3})
\begin{equation}
\|R_0f\|_{\mathcal H(s)}^2\le\|f\|_{\mathcal
H(s)}^2-|f(0)|^2,\quad\forall f\in\mathcal H(s). \label{R0R0}
\end{equation}
Here, the space $\mathcal H(s)$ is finite dimensional and $R_0$
invariant. Thus $R_0$ has a right eigenvector $f$ with eigenvalue
$\overline{a}$; see \cite[p. 36]{MR97h:15020}. Any eigenvector of
$R_0$ is of the form \eqref{sevres-babylone}, and equation
\eqref{R0R0} implies that
\begin{equation}
\label{contrac}
\|f\|^2\le \frac{|f(0)|^2}{1-|a|^2}.
\end{equation}
We will see at the end of the proof of Step 8 that equality in fact holds in \eqref{contrac}.\\

STEP 7: {\sl It holds that $s(a)=0$.}\\

From \cite[p. 282-283]{acs1} it follows that the span of $f$ endowed
with the norm $\|f\|^2=\frac{|f(0)|^2}{1-|a|^2}$ is equal to $\mathcal
H(b_a)$, where $b_a$ is a Blaschke factor, see \eqref{eqBlaschke}. From \eqref{contrac} we get that $\mathcal
H(b_a)$ is contractively included in $\mathcal H(s)$ and from
\cite[Lemme 5.1]{acs1} we then have that the kernel
\begin{equation}
\label{place-voltaire} K_s(p,q)-K_{b_a}(p,q)=\sum_{t=0}^\infty
p^t(b_a(p)\overline{b_a(q)}-s(p)\overline{s(q)})\overline{q}^t
\end{equation}
is positive definite in $\mathbb B_1$. But $b_a(a)=0$. Thus,
setting $p=q=a$ in \eqref{place-voltaire} leads to $s(a)=0$.\\

STEP 8: {\sl We can write $s=b_a\star \sigma_1$, where $\sigma_1$ is a Schur
function.}\\

In the argument we make use of the Hardy space $\mathbf H_2(\mathbb
B_1)$ which is the reproducing kernel Hilbert space with
reproducing kernel
\[
(1-p\overline{q})^{-\star}=\sum_{t=0}^\infty p^t\overline{q}^t.
\]
Note that this is the kernel $k_s$ with $s(p)\equiv 0$.
For more information on this space we refer to \cite{2013arXiv1308.2658A,acs1}.\\

Since a Schur function is bounded in modulus and thus belongs to
the space $\mathbf H_2(\mathbb B_1)$ (see
\cite{2013arXiv1308.2658A}), the representation $s=b_a\star \sigma_1$
with $\sigma_1\in\mathbf H_2(\mathbb B_1)$, follows from \cite[Proof of
Theorem 6.2, p. 109]{MR3127378}. To see that $\sigma_1$ is a Schur
multiplier we note that
\begin{equation}
\label{decom}
K_s(p,q)-K_{b_a}(p,q)=b_a(p)\star K_{\sigma_1}(p,q)\star_r
\overline{b_a(q)}
\end{equation}
implies that $b_a(p)\star K_{\sigma_1}(p,q)\star_r \overline{b_a(q)}$
is positive definite in $\mathbb{B}_1$ and hence $K_{\sigma_1}(p,q)$ is as well
by \cite[Proposition 5.3]{acs3}.\\

STEP 9: {\sl It holds that ${\rm dim}\,(\mathcal H(\sigma_1))={\rm
dim}\,(\mathcal H(s))-1$.}\\

The decomposition \eqref{decom} gives the decomposition
\[
K_s(p,q)=K_{b_a}(p,q)+b_a(p)\star K_{\sigma_1}(p,q)\star_r \overline{b_a(q)}.
\]
The corresponding reproducing kernel spaces do not intersect. Indeed, all elements in the reproducing kernel Hilbert space with reproducing kernel $b_a(p)\star K_{\sigma_1}(p,q)\star_r \overline{b_a(q)}$ vanish at the point $a$ while non zero elements in $\mathcal H(b_a)$ do not vanish. So the decomposition is orthogonal in $\mathcal H(s)$ by Theorem \ref{new}, and equality holds in
\eqref{contrac}. The claim on the dimensions follow.\\

After a finite number of iterations, this procedure leads to a
constant $\sigma_{\ell}$, for some positive integer $\ell$. This constant has to be unitary since the corresponding space
$\mathcal H(\sigma_\ell)$ reduces to $\left\{0\right\}$,  thus proving the theorem.
\end{proof}

We conclude with two remarks and a corollary.

\begin{Rk}
{\rm Given a Blaschke factor the operator of multiplication by $b_a$ is an isometry from ${\bf H}_2(\mathbb{B}_1)$ into itself
(see \cite[Theorem 5.17, p. 106]{MR3127378}), and so is the operator of multiplication by a finite Blaschke product $B$.
The degree of the Blaschke product is the dimension of the space
$\mathbf H_2(\mathbb B_1)\ominus B \mathbf H_2(\mathbb B_1)$.
Thus the previous argument shows in fact that
$\mathcal H(s)$ is isometrically included inside $\mathbf
H_2(\mathbb B_1)$ and that $\mathcal H(s)=\mathbf H_2(\mathbb
B_1)\ominus M_s \mathbf H_2(\mathbb B_1)$.}
\end{Rk}

One can plug a unitary constant $e$ also in the linear fractional
transformation \eqref{lft1} and the same arguments lead to:

\begin{Cy}
If Problem \ref{pb1} has a solution, it is a Blaschke product of
degree ${\rm rank}\, P$.
\end{Cy}

\begin{Rk}
{\rm
The arguments in Steps 5-7 take only into account the fact that the space
$\mathcal H(\Theta)$ is finite dimensional and that $e$ is a unitary constant.
In particular, they also apply in the setting of \cite{2013arXiv1308.2658A},
and in that paper too, the solution of the interpolation problem is a Blaschke
product of degree ${\rm rank}\, P$ when the Pick matrix is degenerate.
}
\end{Rk}

\section{An analogue of Carath\'eodory's theorem in the quaternionic setting}
\label{remarks} Recall first that Carath\'eodory's theorem states
the following (see for instance \cite[pp. 203-205]{burckel},
\cite[p. 48]{sarason94}). We write the result for a radial limit,
but the result holds in fact for a non tangential limit.

\begin{Tm}
Let $s(z)$ be a Schur function and let $e^{it_0}$ be a point on
the unit circle such that
\[
\liminf_{\substack{r\rightarrow
1\\r\in(0,1)}}\frac{1-|s(re^{it_0})|}{1-r}<\infty.
\]
Then, the limits
\[
c=\lim_{\substack{r\rightarrow 1\\r\in(0,1)}}s(re^{it_0})\quad
and\quad \lim_{\substack{r\rightarrow
1\\r\in(0,1)}}\frac{1-s(re^{it_0})\overline{c}}{1-r}
\]
exist, and the second one is positive.
\end{Tm}

This result plays an important role in the classical boundary
interpolation problem for Schur functions. See for instance
\cite{ADLW}, \cite{MR98m:47017}.\\

We prove a related result in the setting of slice-hyperholomorphic
functions. The condition \eqref{place-de-l-opera} will hold
particular for rational functions $s$, as is proved using a
realization of $s$ (see \cite{acs1} for the latter).

\begin{Tm}
Let $s$ be a slice hyperholomorphic Schur function, and assume
that at some point $p_u$ of modulus $1$ we have
\begin{equation}
\label{place-de-la-republique}
\sup_{r\in(0,1)}\frac{1-|s(rp_u)|^2}{1-r^2}<\infty .
\end{equation}
Assume moreover that the function $r\mapsto s(rp_u)$ has a
development in series with respect to the {\sl real} variable $r$
at $r=1$:
\begin{equation}
s(rp_u)=s_u+(r-1)a_u+O(r-1)^2. \label{place-de-l-opera}
\end{equation}
Then
\[
\lim_{\substack{r\rightarrow 1\\ r\in(0,1)}}
 \sum_{t=0}^\infty
r^tp_u^t(1-s(r_np_u)\overline{s_u})\overline{p_u}^t=
(a_u\overline{s_u}-\overline{p_u}a_u\overline{s_u}\, \overline{p_u})(1-\overline{p_u}^2)^{-1}.
\]
\end{Tm}

\begin{proof}
In view of \eqref{place-de-la-republique}, the family of functions $K_s(\cdot, rp_u)$ has a weakly
convergent subsequence. Since weak convergence implies pointwise
convergence the weak limit is readily seen to be the function
$g_u$. Thus
\[
0\le  \langle g_u,g_u \rangle_{\mathcal H(s)}=\lim_{n\rightarrow\infty}
\langle g_u,K_s(\cdot,r_np_u) \rangle_{\mathcal
H(s)}=\lim_{n\rightarrow\infty}g_u(r_np_u),
\]
where $(r_n)_{n\in\mathbb N}$ is a sequence of numbers in $(0,1)$
with limit equal to $1$. Hence we have that
\[
\lim_{n\rightarrow\infty}\sum_{t=0}^\infty r_n^tp_u^t(1-s(r_np_u)\overline{s_u})\overline{p_u}^t\ge 0.
\]
Using \eqref{place-de-l-opera} we have:
\[
\begin{split}
\sum_{t=0}^\infty
r^tp_u^t(1-s(r_np_u)\overline{s_u})\overline{p_u}^t&=
\sum_{t=0}^\infty r^tp_u^t((r-1)a_u\overline{s_u}+O(r-1)^2)\overline{p_u}^t\\
&=((r-1)a_u\overline{s_u}-r\overline{p_u}(r-1)a_u\overline{s_u}\,
\overline{p_u})(1-r)^{-1}(1-r\overline{p_u}^2)^{-1}+\\
&\hspace{5mm}+
\sum_{t=0}^\infty r^tp_u^nO(r-1)^2\overline{p_u}^t\\
&=(a_u\overline{s_u}-r\overline{p_u}a_u\overline{s_u}\, \overline{p_u})
(1-r\overline{p_u}^2)^{-1}+\\
&\hspace{5mm}+ \sum_{t=0}^\infty r^tp_u^nO(r-1)^2\overline{p_u}^t.
\end{split}
\]
This expression tends to
\begin{equation}
\label{barbes}
(a_u\overline{s_u}-\overline{p_u}a_u\overline{s_u}\, \overline{p_u})(1-\overline{p_u}^2)^{-1},
\end{equation}
as $r\rightarrow 1$.
\end{proof}
\begin{Rk}{\rm
The example $s(p)=\frac{1+pa}{2}$, where $a\in\mathbb B_1$ is such
that $ap_u\not=p_ua$, shows that \eqref{barbes} is different, in
general, from $a_u\overline{s_u}$.}
\end{Rk}
\bibliographystyle{plain}

\begin{thebibliography}{10}

\bibitem{2013arXiv1308.2658A}
D.~{Alpay}, V.~{Bolotnikov}, F.~{Colombo}, and I.~{Sabadini}.
\newblock {Self-mappings of the quaternionic unit ball: multiplier properties,
  Schwarz-Pick inequality, and Nevanlinna--Pick interpolation problem}.
\newblock {\em ArXiv e-prints}, August 2013.
\newblock To appear in the Indiana Mathematical Journal of Mathematics.

\bibitem{abds2}
D.~Alpay, P.~Bruinsma, A.~Dijksma, and {H.S.V. de} Snoo.
\newblock Interpolation problems, extensions of symmetric operators and
  reproducing kernel spaces {II}.
\newblock {\em Integral Equations Operator Theory}, 14:465--500, 1991.

\bibitem{abds3}
D.~Alpay, P.~Bruinsma, A.~Dijksma, and {H.S.V. de} Snoo.
\newblock Interpolation problems, extensions of symmetric operators and
  reproducing kernel spaces {II} (missing section 3).
\newblock {\em Integral Equations Operator Theory}, 15:378--388, 1992.

\bibitem{2013arXiv1310.1035A}
D.~{Alpay}, F.~{Colombo}, I.~{Lewkowicz}, and I.~{Sabadini}.
\newblock {Realizations of slice hyperholomorphic generalized contractive and
  positive functions}.
\newblock {\em ArXiv e-prints}, October 2013.

\bibitem{acs3}
D.~{Alpay}, F.~{Colombo}, and I.~{Sabadini}.
\newblock {Krein-Langer factorization and related topics in the slice
  hyperholomorphic setting}.
\newblock Journal of Geometric Analysis. Accepted. To appear.

\bibitem{acs1}
D.~{Alpay}, F.~{Colombo}, and I.~{Sabadini}.
\newblock {Schur functions and their realizations in the slice hyperholomorphic
  setting}.
\newblock {\em Integral Equations and Operator Theory}, 72:253--289, 2012.

\bibitem{MR3127378}
D.~Alpay, F.~Colombo, and I.~Sabadini.
\newblock Pontryagin-de {B}ranges-{R}ovnyak spaces of slice hyperholomorphic
  functions.
\newblock {\em J. Anal. Math.}, 121:87--125, 2013.

\bibitem{ADLW}
D.~Alpay, A.~Dijksma, H.~Langer, and G.~Wanjala.
\newblock {Basic boundary interpolation for generalized {S}chur functions and
  factorization of rational $J$--unitary matrix functions}.
\newblock In D.~Alpay and I.~Gohberg, editors, {\em {Interpolation, Schur
  functions and moment problems}}, volume 165 of {\em Oper. Theory Adv. Appl.},
  pages 1--29. Birkh{\" a}user Verlag, Basel, 2006.

\bibitem{adubi1}
D.~Alpay and C.~Dubi.
\newblock Boundary interpolation in the ball.
\newblock {\em Linear {A}lgebra and {A}pplications}, 340:33--54, 2002.

\bibitem{ad-laa-herm}
D.~Alpay and H.~Dym.
\newblock {Structured invariant spaces of vector valued functions, hermitian
  forms, and a generalization of the Iohvidov laws}.
\newblock {\em Linear Algebra Appl.}, 137/138:137--181, 1990.

\bibitem{ad-jfa}
D.~Alpay and H.~Dym.
\newblock On a new class of structured reproducing kernel {H}ilbert spaces.
\newblock {\em J. Funct. Anal.}, 111:1--28, 1993.

\bibitem{aron}
N.~Aronszajn.
\newblock Theory of reproducing kernels.
\newblock {\em Trans. Amer. Math. Soc.}, 68:337--404, 1950.

\bibitem{ball-contrac}
J.~Ball.
\newblock Models for noncontractions.
\newblock {\em J. Math. Anal. Appl.}, 52:235--259, 1975.

\bibitem{dbhsaf1}
{L. de} Branges.
\newblock Some {Hilbert} spaces of analytic functions {I}.
\newblock {\em Trans. Amer. Math. Soc.}, 106:445--468, 1963.

\bibitem{burckel}
R.B. Burckel.
\newblock {\em An introduction to classical complex analysis, {V}ol. 1}.
\newblock Birkh{\"a}user, 1979.

\bibitem{MR2752913}
F. Colombo, I. Sabadini, and D.~C. Struppa.
\newblock {\em Noncommutative functional calculus}, volume 289 of {\em Progress
  in Mathematics}.
\newblock Birkh\"auser/Springer Basel AG, Basel, 2011.
\newblock Theory and applications of slice hyperholomorphic functions.

\bibitem{Dym_CBMS}
H.~Dym.
\newblock {\em ${J}$--contractive matrix functions, reproducing kernel
  {H}ilbert spaces and interpolation}.
\newblock Published for the Conference Board of the Mathematical Sciences,
  Washington, DC, 1989.

\bibitem{fw}
P.A. Fillmore and J.P. Williams.
\newblock On operator ranges.
\newblock {\em {Advances in Mathematics}}, 7:254--281, 1971.

\bibitem{HornJohnson}
R.~A. Horn and C.~R. Johnson.
\newblock {\em Topics in matrix analysis}.
\newblock Cambridge University Press, Cambridge, 1994.
\newblock Corrected reprint of the 1991 original.

\bibitem{kky}
V.~Katsnelson, A.~Kheifets, and P.~Yuditskii.
\newblock An abstract interpolation problem and the extension theory of
  isometric operators.
\newblock In H.~Dym, B.~Fritzsche, V.~Katsnelson, and B.~Kirstein, editors,
  {\em Topics in interpolation theory}, volume~95 of {\em {Operator {T}heory:
  {A}dvances and {A}pplications}}, pages 283--297. Birkh{\" a}user Verlag,
  Basel, 1997.
\newblock Translated from: Operators in function spaces and problems in
  function theory, p. 83--96 ({N}aukova--{D}umka, {K}iev, 1987. {E}dited by
  {V}.{A}. {M}archenko).

\bibitem{HM}
J.~Rovnyak.
\newblock Characterization of spaces ${H(M)}$.
\newblock Unpublished paper, $1968$. Available at the URL {\tt
  http://www.people.virginia.edu/~jlr5m/home.html}.

\bibitem{sarason94}
D.~Sarason.
\newblock {\em Sub--{H}ardy {H}ilbert spaces in the unit disk}, volume~10 of
  {\em University of {A}rkansas lecture notes in the mathematical sciences}.
\newblock Wiley, {N}ew {Y}ork, 1994.

\bibitem{MR98m:47017}
D.~Sarason.
\newblock Nevanlinna-{P}ick interpolation with boundary data.
\newblock {\em Integral Equations Operator Theory}, 30(2):231--250, 1998.
\newblock Dedicated to the memory of Mark Grigorievich Krein (1907--1989).

\bibitem{schwartz}
L.~Schwartz.
\newblock Sous espaces hilbertiens d'espaces vectoriels topologiques et noyaux
  associ\'{e}s (noyaux reproduisants).
\newblock {\em J. Analyse Math.}, 13:115--256, 1964.

\bibitem{MR97h:15020}
F.~Zhang.
\newblock Quaternions and matrices of quaternions.
\newblock {\em Linear Algebra Appl.}, 251:21--57, 1997.

\end{thebibliography}
\def\cprime{$'$} \def\lfhook#1{\setbox0=\hbox{#1}{\ooalign{\hidewidth
  \lower1.5ex\hbox{'}\hidewidth\crcr\unhbox0}}} \def\cprime{$'$}
  \def\cfgrv#1{\ifmmode\setbox7\hbox{$\accent"5E#1$}\else
  \setbox7\hbox{\accent"5E#1}\penalty 10000\relax\fi\raise 1\ht7
  \hbox{\lower1.05ex\hbox to 1\wd7{\hss\accent"12\hss}}\penalty 10000
  \hskip-1\wd7\penalty 10000\box7} \def\cprime{$'$} \def\cprime{$'$}
  \def\cprime{$'$} \def\cprime{$'$} \def\cprime{$'$}

\end{document}